\documentclass[11pt, english]{smfart}
\usepackage{amsfonts, amssymb, amsmath, latexsym}
\usepackage{mathptmx}
\usepackage[english]{babel}
\usepackage{mysectionsen}
\usepackage[latin1]{inputenc}
\usepackage[T1]{fontenc}
\usepackage{textcomp}
\usepackage[dvips]{graphicx}
\usepackage{letltxmacro}

\newtheorem{Conj}{Conjecture}

\input xy
\xyoption{all}

\begin{document}

\thispagestyle{empty}

\begin{center}
{\Large \bfseries{\textsc{
On some pro-$p$ groups \\
\vspace{0.2cm}
from infinite-dimensional Lie theory}}}
\end{center}

\vspace{1cm}

\begin{center}
\textsc{Inna Capdeboscq and Bertrand R\'emy}
\end{center}

\vspace{1cm}

\begin{center} \today
\end{center}

\vspace{5cm}

\hrule

\vspace{0.5cm}

{\small
\noindent
{\bf Abstract:}
We initiate the study of some pro-$p$-groups arising from infinite-dimensional Lie theory.
These groups are completions of some subgroups of incomplete Kac-Moody groups over finite fields, with respect to various completions of algebraic or geometric origin. 
We show topological finite generation for the pro-$p$ Sylow subgroups in many complete Kac-Moody groups.
This implies abstract simplicity of the latter groups.
We also discuss with the question of (non-)linearity of these pro-$p$ groups.

\vspace{0,2cm}


\vspace{0,1cm}

\noindent
{\bf AMS classification (2000):}
17B22,
17B67,
20E32,
20E42,
20F20, 
20G44,
51E24.
}

\vspace{0,5cm}

\hrule

\newpage

\vspace{1cm}

\tableofcontents

\newpage

\section*{Introduction}

The general theme of this paper is a connection between infinite-dimensional Lie theory and pro-$p$ groups.
More precisely, we are interested in Kac-Moody groups over finite fields: in various complete versions of these groups, we obtain locally compact, totally disconnected groups admitting a useful action on an explicit building.
This action and some metric properties of the building show that if the ground field is finite, say of characteristic $p$, then maximal pro-$p$ subgroups exist and are conjugate in these (non-compact) groups -- this is why we allow ourselves to call these subgroups the {\it pro-$p$ Sylow subgroups}~of the complete Kac-Moody group under consideration.

\medskip 

This situation is a generalization of the action of the non-Archimedean simple Lie group ${\rm SL}_n\bigl( \mathbf{F}_q (\!(t)\!) \bigr)$ on its Bruhat-Tits building \cite{Corvallis}.
In this case, the interesting pro-$p$ Sylow subgroup is, up to conjugation, the subgroup of the matrices of ${\rm SL}_n( \mathbf{F}_q [[t]])$ which are unipotent upper triangular modulo $t$.
In the general case, the ambient locally compact groups we consider are non-linear, at least because they often contain infinite finitely generated simple subgroups.
Nevertheless, many challenging questions arise from the analogy with the previous classical situation.

\medskip

In this paper, we address the question of topological finite generation of the pro-$p$ Sylow subgroups and the question of abstract simplicity for the complete Kac-Moody groups.
Here are simplified versions of the results we obtain.
Recall that a Kac-Moody group is basically defined by a generalized Cartan matrix and a ground field (see Sect. \ref{s - AlgComp}).

\medskip

\noindent
\emph{\textbf{Theorem 1}} --- \emph{When the characteristic $p$ of the finite ground field is greater than the opposite of any off-diagonal coefficient of the generalized Cartan matrix, the pro-$p$ Sylow subgroups of complete Kac-Moody groups are topologically finitely generated; namely, they are generated by the simple root groups.}

\medskip

For a more precise version, see Theorem \ref{th.Frattini} and Corollary \ref{cor.FGRR} where we actually compute the Frattini subgroup, hence the first homology group with coefficients in ${\bf Z}/p{\bf Z}$ , of the above pro-$p$ groups.
Similar finite generation cases, with restrictions on the type of generalized Cartan matrix, had already been obtained in \cite{AbrMuh} and \cite{CER}. 

\medskip

\noindent
\emph{\textbf{Theorem 2}} --- \emph{Under the same conditions as in the previous statement, natural subquotients of complete Kac-Moody groups over finite fields are simple as abstract groups.}

\medskip

For a more precise version, see Theorem \ref{th.simple}.
We note that T.~Marquis obtained a better result since he needs no assumption on the size of the ground field \cite[Theorem A]{Marquis}; his techniques are of dynamical nature (contraction groups).
Our proof is inspired by Carbone-Ershov-Ritter's proof of abstract simplicity for an already wide family of complete Kac-Moody groups over finite fields \cite{CER}; it is algebraic in nature since it is a combination of techniques from Tits systems and pro-$p$ groups (a combination previously used, but not in its full strength, for the proof of topological simplicity in \cite{RemyGAFA}).

\medskip

Concerning the question of linearity, this is the approach we chose to support the fact that the pro-$p$ groups we are starting to investigate in this paper are new.
The results we obtain are not optimal, and we hope to improve them in a near future by using deeper techniques from pro-$p$ groups.

\medskip

Mentioning this project is the opportunity to mention also a further motivation to study pro-$p$ Sylow subgroups in locally compact Kac-Moody groups.
It is now well-known that some Lie-theoretic techniques are very useful in the general study of pro-$p$ groups \cite{horizons}.
The underlying (not yet completely exploited) idea of this paper is to show that conversely, (infinite-dimensional) Lie theory provides a lot of interesting pro-$p$ groups which perfectly fit in this approach of studying pro-$p$ groups by making emerge algebraic structures; indeed, these algebraic structures already exist here by construction!

\medskip

Finally, we note that Kac-Moody groups have already provided interesting pro-$p$ groups: for instance, M.~Ershov proved the Golod-Shafarevich property for the full pro-$p$ completion of some finitely generated (actually Kazhdan) groups arising as finite index subgroups of Borel subgroups in suitable minimal Kac-Moody groups over finite fields \cite{ErshovGoSha}.
Here, we investigate pro-$p$ groups which are not full pro-$p$ completions, but which come from intermediate completion procedures, both of algebraic \cite{Rousseau} and of geometric nature \cite{RemRon}. 

\medskip

\textbf{Structure of the paper.}
Section 1 is dedicated to the presentation of Kac-Moody objects: we start from Lie algebras, then we sum very quickly J.~Tits' construction of the minimal version of Kac-Moody, and we finally introduce G.~Rousseau's very recent approach to "algebraic completions" of the previous groups.
Section 2 discusses the problem of topological finite generation for the pro-$p$ Sylow subgroups of the latter groups over finite fields of characteristic $p$.
Section 3 exploits the previous finite generation results to derive some abstract simplicity statements for various versions of complete Kac-Moody groups.
At last, Section 4 initiates the problem of linearity of the pro-$p$ groups we are interested in.
Each of the last three sections contains a last subsection mentioning some open problems in the field.

\medskip

\textbf{Conventions.}
In this paper we adopt the following notation.
If $G$ is any group (topological or not), its group of commutator is denoted by $[G,G]$ and the characteristic subgroup of all the $k$-th powers of all elements in $G$ is denoted by $G^k$.
By a {\it local field}~we mean a field endowed with a discrete valuation and complete for the corresponding absolute value.

\vspace{1cm}

\section{Algebraic completions of Kac-Moody groups}
\label{s - AlgComp}

In this expository section, we quickly introduce some objects from Kac-Moody theory, from Lie algebras and root systems to more and more sophisticated group functors.
As a general reference where this material is also presented in a comprehensive paper, we recommend J. Tits' Bourbaki seminar talk \cite{TitsBBK}. 

\subsection{Kac-Moody Lie algebras}
\label{ss - KM Lie}

The reference for this section is \cite{TitsKM}.

\medskip 

Let us start with the Lie-theoretic data needed to define the objects we are interested in. 
The basic definition is that of a  {\it generalized Cartan matrix}: it is an integral matrix $A=[A_{s,t}]_{s,t \in S}$ satisfying: $A_{s,s}=2$, $A_{s,t} \leqslant  0$ when $s \neq t$ and $A_{s,t}=0 \Leftrightarrow A_{t,s}=0$. 
A more precise data is given by a {\it Kac-Moody root datum}, namely a $5$-tuple  $\mathcal{D}=\bigl(S, A, \Lambda, (c_s)_{s \in S}, (h_s)_{s \in S} \bigr)$, where $A$ is a  generalized Cartan matrix indexed by the finite set $S$, $\Lambda$ is a free ${\bf Z}$-module (whose ${\bf Z}$-dual is denoted by $\Lambda^\vee$); the elements $c_s$ of $\Lambda$ and $h_s$ of $\Lambda^\vee$ are required to satisfy $c_s(h_t) =A_{ts}$ for all $s$ and $t$ in $S$. 
From such a $\mathcal{D}$, one defines a complex Lie algebra, say $\frak{g}_\mathcal{D}$,  thanks to a presentation formally generalizing Serre's presentation of finite-dimensional semisimple Lie algebras.
This presentation involves canonical generators $e_s$ and $f_s$ for $s \in S$.

\medskip

Next we introduce the free abelian group $Q = \bigoplus_{s \in S} {\bf Z}Ê\alpha_s$ on the symbols $\alpha_s$ (with same index set as for $A$).
We can then define on the Lie algebra $\frak{g}_\mathcal{D}$ a $Q$-gradation in which the degrees with non-trivial corresponding spaces belong to $Q^+ \cup Q^-$, where $Q^+ = \sum_{s \in S} {\bf N}Ê\alpha_s$ and $Q^-=-Q^+$.
The latter non-zero degrees are called {\it roots}; when the family $(c_s)_{s \in S}$ is free over ${\bf Z}$ in $\Lambda$, the defining relations of $\frak{g}_\mathcal{D}$ lead to an interpretation of roots in terms of the restriction to a maximal abelian subalgebra of the adjoint action of $\frak{g}_\mathcal{D}$ on itself.
The {\it height} of a root $\alpha = \sum_{s \in S} n_s \alpha_s$ is the integer ${\rm ht}(\alpha) = \sum_{s \in S} n_s$. 
Moreover there is a natural ${\bf Z}$-linear action on the lattice $Q$ by a Coxeter group $W$ generated by involutions denoted by $s \in S$ (same index set as before): the action is defined by the formula $s.a_t=a_t - A_{st}a_s$.
Roots are permuted by this action, and a root is called {\it real} if it is in the $W$-orbit of a {\it simple} root, i.e. some $\alpha_s$; a root is called {\it imaginary}~otherwise.
The set of roots (resp. real roots, imaginary roots) is denoted by $\Delta$ (resp. $\Delta_{\rm re}$, $\Delta_{\rm im}$).

\medskip

As usual, one  attaches to $\frak{g}_\mathcal{D}$ its universal envelopping algebra \cite[Part I, Chapter III]{SerreLie}, denoted by $\mathcal{U}\frak{g}_\mathcal{D}$.
Using (among other things) divided powers ${1 \over n!}Êe_s^n$ and ${1 \over n!}Êf_s^n$ of the canonical generators $e_s$ and $f_s$ and of their Weyl group conjugates, J. Tits defined a certain ${\bf Z}$-form of 
$\mathcal{U}\frak{g}_\mathcal{D}$, denoted by $\mathcal{U}_\mathcal{D}$: it is a $Q$-filtered ring which he used to define the minimal group functors below, and whose completions with respect to some subsemigroups are used in the definition and study of some complete Kac-Moody groups (\ref{ss - KM complete}).

\subsection{Minimal Kac-Moody groups}
\label{ss - KM Tits}

A reference for this subsection may be \cite{RemyAst}.

\medskip 

The adjoint action of $\frak{g}_\mathcal{D}$ on itself extends to an action on $\mathcal{U}\frak{g}_\mathcal{D}$ in which the elements in real root spaces have a locally nilpotent action. 
Therefore the latter actions can be exponentiated and give rise to 1-parameter "unipotent" subgroups in the automorphism group of $\mathcal{U}\frak{g}_\mathcal{D}$, and actually in the 
automorphism group of the ${\bf Z}$-form $\mathcal{U}_\mathcal{D}$ for suitable restrictions of parameters and elements in $\frak{g}_\mathcal{D}$.
Roughly speaking, J. Tits defined a {\it minimal Kac-Moody group functor}~$\frak{G}_\mathcal{D}$ as an amalgamation of a torus defined by $\Lambda$ and a quotient of the subgroup generated by these 1-parameter subgroups.
To each real root $\gamma \in \Delta_{\rm re}$ is attached a subgroup functor $\mathcal{U}_\gamma$, which is isomorphic to the 1-dimension additive group functor.
Nevertheless, there is no subgroup attached to an imaginary root in this (minimal) version of Kac-Moody groups.

\medskip

\begin{Ex}
\label{Ex.minKM}
The most standard example of a minimal Kac-Moody group functor comes from the generalized Cartan matrix of (affine) type $\widetilde{A}_n$: it is the group functor which attaches to any field $k$ the group ${\rm SL}_{n+1}(k[t,t^{-1}])$.
\end{Ex}

\medskip

We will not go into details for these groups in this paper, but we simply recall that for each field $k$, the group $\frak{G}_\mathcal{D}(k)$ has a nice action on pair of two twinned buildings.
These buildings are cell complexes which admit a very useful complete, non-positively curved metric.
For this distance, we have a fixed point lemma for isometric group actions with bounded orbits.
The buildings are locally finite if and only if the ground field $k$ is finite, in which case the full isometry groups of the two buildings are locally compact.
Note that thanks to techniques relevant to discrete subgroups of Lie groups and Coxeter groups, the question of simplicity for the groups $\frak{G}_\mathcal{D}(k)$, where $k$ is finite, is almost completely solved (see \cite{CapRem} and \cite{CapRemBis}). 

\subsection{Complete Kac-Moody groups}
\label{ss - KM complete}

We saw in the previous section that minimal Kac-Moody groups are (non-linear) generalizations of groups of the form ${\bf G}(k[t,t^{-1}])$ where $k$ is a field and ${\bf G}$ is an $F$-isotropic semisimple group over $k$.
We want to present various complete versions of Kac-Moody groups which generalize, in a similar way, groups like ${\bf G}\bigl(k(\!(t)\!)\bigr)$.
Some constructions are due to O.~Mathieu \cite{Mathieu} or Sh.~Kumar \cite{Kumar} over the complex numbers; we present here a refined version due to G.~Rousseau, whose paper \cite{Rousseau} is the reference for this subsection.
He obtains here group functors defined over ${\bf Z}$ carrying a structure of ind-scheme.

\medskip 

The starting point is the fact (from group scheme theory) that for suitable affine groups, the corresponding rings of regular functions can be seen as restricted duals of the universal enveloping algebra of the corresponding Lie algebra.
In order to construct groups over ${\bf Z}$, the idea is therefore to go the other way round and to start from a suitable ${\bf Z}$-form of an enveloping algebra. 
Actually, for a given Kac-Moody  root datum $\mathcal{D}$, the groups that are constructed in this context are obtained from closed sub-rings of the completion of $\mathcal{U}_\mathcal{D}$ with respect to the elements in the positive part of the $Q$-gradation.
For each element $x \in \frak{g}_\mathcal{D} \cap \mathcal{U}_\mathcal{D}$, G.~Rousseau makes a choice of {\it exponential sequence}~$(x^{[n]})_{n \geqslant 1}$ (not unique, but whose non-uniqueness is quite well-controlled) \cite[Propositions 2.4 and 2.7]{Rousseau}.
From this choice an exponential can be associated to $x$ in a way which is compatible both with respect to a natural coproduct and to the $Q$-filtration: this enables Rousseau to associate to any closed set of roots, say $\Psi$, a very useful pro-unipotent group scheme over ${\bf Z}$, denoted by $\frak{U}_\Psi^{\rm ma}$ \cite[3.1]{Rousseau}.
These groups have a nice filtration in terms of the root system, and note that in this context there exist root groups associated to imaginary roots: the fact that imaginary roots are taken into account in this construction is a major tool in our paper.
Note that $k$ is a finite field of characteristic $p$, the group of $k$-points $\frak{U}_\Psi^{\rm ma}(k)$ is naturally a pro-$p$ group.

\medskip

The construction of the full complete Kac-Moody group $\frak{G}_\mathcal{D}^{\rm ma+}$ is not so useful in what follows.
Once again, we are mainly interested in the geometric outcome of the construction, namely the fact that for any field $k$, the group $\frak{G}_\mathcal{D}^{\rm ma+}(k)$ has a natural strongly transitive action on a building whose apartments are described thanks the to Weyl group $W$ of $\mathcal{D}$.
The building is closely related to the twin buildings of the previous subsections \cite[Corollaire 3.18]{Rousseau}. 
In fact this geometric viewpoint provides a connection with another, less subtle, completion defined in \cite[1.B]{RemRon} when the ground field $k$ is finite.
In this case, the completion defined there, and denoted by $\frak{G}_\mathcal{D}^{\rm rr}(k)$, is simply the closure of the image of $\frak{G}_\mathcal{D}(k)$ in the (locally compact) full automorphism group of the positive twin building (\ref{ss - KM Tits}). 

\medskip

\textbf{Notation.}
{\it In the rest of the paper, we are given a Kac-Moody root datum $\mathcal{D}$ and a finite field $\mathbf{F}_q$ of order $q$ and of characteristic $p$.
We are mainly interested in the locally pro-$p$ Kac-Moody group $\frak{G}_\mathcal{D}^{\rm ma+}(\mathbf{F}_q)$, which we also call the {\rm algebraic completion}~of Tits' minimal Kac-Moody group $\frak{G}_\mathcal{D}(\mathbf{F}_q)$.}

\section{Frattini subgroups of  $p$-Sylow subgroups}
\label{s - TopFG}

As already mentioned, the family given by the (usually non-linear) groups $\frak{G}_\mathcal{D}^{\rm ma+}(\mathbf{F}_q)$ can be seen as a vast generalization of the family of groups ${\bf G}\bigl( \mathbf{F}_q(\!(t)\!)\bigr)$ where $\mathbf{G}$ is a semisimple algebraic group.
Moreover the argument given in \cite[1.B.2]{RemyGAFA} implies that the conjugates of the subgroup $\frak{U}_{\Delta^+}^{\rm ma}(\mathbf{F}_q)$ (where $\Delta^+$ is the closed set of all positive roots of $\mathcal{D}$) can be group-theoretically characterized as the pro-$p$ Sylow subgroups of $\frak{G}_\mathcal{D}^{\rm ma+}(\mathbf{F}_q)$.
Therefore, since  when $\frak{G}_\mathcal{D}^{\rm ma+}(\mathbf{F}_q) = {\rm SL}_n\bigl( \mathbf{F}_q(\!(t)\!)\bigr)$ the group $\frak{U}_{\Delta^+}^{\rm ma}(\mathbf{F}_q)$ has finite index in ${\rm SL}_n(\mathbf{F}_q[[t]])$, it is natural to ask whether an arbitrary group $\frak{U}_{\Delta^+}^{\rm ma}(\mathbf{F}_q)$ shares some interesting properties with the more classical pro-$p$ Sylow subgroups of the groups ${\rm SL}_n(\mathbf{F}_q[[t]])$.
We focus here on topological finite generation, via the computation of the Frattini subgroups of $\frak{U}_{\Delta^+}^{\rm ma}(\mathbf{F}_q)$.

\medskip

\textbf{Notation.}
{\it For the sake of simplicity, we use the notation $U^{{\rm ma}+}$ instead of $\frak{U}_{\Delta^+}^{\rm ma}(\mathbf{F}_q)$ in this section.}

\subsection{Some facts from pro-$p$ groups}
\label{ss - pro-p basics}

We collect here some facts about pro-$p$ groups that we will freely use in the next subsection.

\medskip

Let $V$ be a profinite group.
Then its {\it Frattini subgroup}, denoted by $\Phi(V)$, is by definition the intersection of all its maximal proper open subgroups.
It is a closed normal subgroup in $V$, and it can be proved that a subset $X$ in $V$ topologically generates $V$ if and only if its image in $V/\Phi(V)$ generates the latter group topologically
\cite[Proposition 1.9]{DDMS}.

\medskip

We henceforth assume that $V$ is pro-$p$.

\medskip

This implies in particular that $\Phi(V) = \overline{V^p[V,V]}$, where $V^p$ denotes the characteristic subgroup of all $p$-th powers in $V$ \cite[Proposition 1.13]{DDMS}.
It is also well-known that there is a natural isomorphism ${\rm H}_1(V,{\bf Z}/p{\bf Z}) \cong V/\Phi(V)$ \cite[Lemma 6.8.6]{RZ}.

\subsection{Finite generation}
\label{ss - fg}
We prove here a refined version of topological finite generation for maximal pro-$p$ subgroups, namely we provide a Lie-theoretic description of the relevant Frattini subgroups.

\smallskip

We start by recalling a result which already appeared in the literature as \cite[Lemma 7.1]{CER}.

\medskip

\begin{Lemma}
\label{lemma.CER}
Let $P$ be a topological group.
Let $(K_i)_{i \geqslant 1}$ be a sequence of open normal subgroups which is a basis of neighborhoods of the identity element.
Let $V$ be a compact subgroup that satisfies the inclusion 

\smallskip

\centerline{$(*_i)$ \quad $K_i \subseteq V \cdot K_{i+1}$.}

\smallskip

\noindent for each $i \geqslant 1$.
Then we have $K_1 \subseteq V$. 
\end{Lemma}

\medskip

\noindent \emph{\textbf{Proof}}.
Let $g_1 \in K_1$. 
We can successively write: $g_1=v_1'g_2$ for $v_1' \in V$ and $g_2 \in K_2$, $g_2=v_2'g_3$ for $v_2' \in V$ and $g_3 \in K_3$, $\dots$ $g_i=v_i'g_{i+1}$ for $v_i' \in V$ and $g_{i+1} \in K_{i+1}$, $\dots$ so that setting $v_i=v_1'v_2'\cdots v_i'$ we obtain a sequence of elements $(v_i)_{i \geqslant 1}$ in $V$ such that $g_1 = v_ig_{i+1}$ for each $i \geqslant 1$. 
By compactness of $V$ and since $(K_i)_{i \geqslant 1}$ is a basis of neighborhoods of the identity element, we deduce that $g_1 \in V$.
\hfill $\Box$

\medskip

We can now turn to the description of the Frattini subgroup of some maximal pro-$p$ subgroups in completed Kac-Moody groups.

\medskip 

\begin{Thm}
\label{th.Frattini}
Let $\frak{G}_\mathcal{D}^{\rm ma+}({\bf F}_q)$ be the algebraic completion of the minimal Kac-Moody group $\frak{G}_\mathcal{D}({\bf F}_q)$.
Let $A$ be the generalized Cartan matrix involved in the Kac-Moody root datum $\mathcal{D}$ defining the group functor $\frak{G}_\mathcal{D}$.
Let $U^{{\rm ma}+}$ be the standard pro-$p$ Sylow subgroup of $\frak{G}_\mathcal{D}^{\rm ma+}({\bf F}_q)$.
We assume that the characteristic $p$ of ${\bf F}_q$ is greater than the absolute value of any off-diagonal coefficient of $A$.
Then the following holds.

\begin{itemize}
\item[{\rm (i)}] The Frattini subgroup $\Phi(U^{{\rm ma}+})$ is equal to the derived group $[U^{{\rm ma}+},U^{{\rm ma}+}]$.
\item[{\rm (ii)}] We have: $\Phi(U^{{\rm ma}+}) = \overline{\langle U_\gamma : \gamma  \hbox{\rm \, non-simple positive real root} \rangle}$.
\item[{\rm (iii)}] We have: ${\rm H}_1(U^{{\rm ma}+},{\bf Z}/p{\bf Z}) \simeq ({\bf Z}/p{\bf Z})^{{\rm size}(A) \cdot [{\bf F}_q:{\bf Z}/p{\bf Z}]}$.
\end{itemize}

\noindent In particular, the group $U^{{\rm ma}+}$ is a topologically finitely generated pro-$p$ group.
\end{Thm}

\medskip 

\noindent Note that by \cite[Proposition 6.11]{Rousseau} and under the assumptions of the theorem, the group $U^{{\rm ma}+}$ is the closure in $\frak{G}_\mathcal{D}^{\rm ma+}({\bf F}_q)$ of the subgroup of $\frak{G}_\mathcal{D}({\bf F}_q)$ generated by the root groups indexed by the positive real roots.

\medskip

\begin{Rk}
\label{rk.generating.set}
We denote by $\{\alpha_s \}_{s \in S}$ the simple roots of the Kac-Moody root datum $\mathcal{D}$ and by $u_\alpha : ({\frak{g}_\alpha}_{{\bf F}_q} \simeq) \, {\bf F}_q \to U_\alpha$ the additive parametrization of each root group in Tits' presentation. 
We assume that the ground field ${\bf F}_q$ is of degree $r$ over ${\bf Z}/p{\bf Z}$ and pick a corresponding basis $v_1=1$, $v_2$, $\dots$ $v_r$.
Then a topologically generating set is concretely given by the elements $u_{\alpha_s}(v_i)$ with $s \in S$ and $1 \leqslant i \leqslant r$. 
This (fortunately!) agrees with the interpretation of ${\rm H}_1(U^{{\rm ma}+},{\bf Z}/p{\bf Z})$ in terms of topologically generating sets  \cite[I.4.2]{SerreCoGa}.
\end{Rk}

\medskip

\noindent \emph{\textbf{Proof}}.
This is precisely the proof of \cite[Proposition 6.11]{Rousseau} that we combine with the previous lemma to obtain the theorem.
More precisely, using Rousseau's notation, let us set $V = \overline{[U^{{\rm ma}+},U^{{\rm ma}+}]}$ and $K_i=U^{{\rm ma}+}_{i+1}$.
What is checked in the proof is precisely that $(*_i)$ is satisfied for each $i \geqslant 1$.
The lemma hence implies that $\overline{[U^{{\rm ma}+},U^{{\rm ma}+}]}$ contains $U^{{\rm ma}+}_2$.
But by \cite[3.4]{Rousseau}, we know that the group $U^{{\rm ma}+}/U^{{\rm ma}+}_2$ is isomorphic to the finite-dimensional commutative Lie algebra over ${\bf F}_q$ equal to $\bigoplus_{\gamma {\rm \,\, simple}} {\frak{g}_\gamma}_{{\bf F}_q}$.
Since the latter group is an elementary $p$-group, we conclude that we have the equalities $\Phi(U^{{\rm ma}+}) = \overline{[U^{{\rm ma}+},U^{{\rm ma}+}]} = U^{{\rm ma}+}_2$.
This finally implies the theorem, since we eventually have $\overline{[U^{{\rm ma}+},U^{{\rm ma}+}]}=[U^{{\rm ma}+},U^{{\rm ma}+}]$ by \cite[Proposition 1.19]{DDMS}.
\hfill $\Box$

\medskip

\begin{Cor}
\label{cor.FGRR}
Let $\overline{U}^{\rm rr}$ be the standard pro-$p$ Sylow subgroup in the geometric completion of the minimal Kac-Moody group $\frak{G}_\mathcal{D}({\bf F}_q)$.
Let $\overline{U}^{{\rm cg}\lambda}$ be the standard pro-$p$ Sylow subgroup in the representation-theoretic completion of the minimal Kac-Moody group $\frak{G}_\mathcal{D}({\bf F}_q)$ associated to the regular dominant weight $\lambda$.
Then, under the same hypothesis on $p$ and $A$ as in the theorem, the same statements hold for each of the two pro-$p$ groups $\overline{U}^{\rm rr}$ and $\overline{U}^{{\rm cg}\lambda}$ replacing $U^{{\rm ma}+}$.
\end{Cor}

\medskip

\noindent \emph{\textbf{Proof}}.
By \cite[6.3 5) and Proposition 6.11]{Rousseau}, there are continuous, proper, open surjective homomorphisms $U^{{\rm ma}+} \twoheadrightarrow \overline{U}^{{\rm cg}\lambda}\twoheadrightarrow \overline{U}^{\rm rr}$.
We denote by $\varphi$ the composed map.
First, in view of the behaviour of Frattini subgroups under epimorphisms of pro-$p$ groups \cite[Corollary 2.8.8]{RZ}, we obtain (i) and (ii) in the above theorem for $\overline{U}^{\rm rr}$ and $\overline{U}^{{\rm cg}\lambda}$ replacing $U^{{\rm ma}+}$.
The same reference also implies that $\varphi$ induces an epimorphism ${\rm H}_1(U^{{\rm ma}+},{\bf Z}/p{\bf Z}) \twoheadrightarrow {\rm H}_1(\overline{U}^{\rm rr},{\bf Z}/p{\bf Z})$. 
Therefore the dimension over ${\bf Z}/p{\bf Z}$ of ${\rm H}_1(\overline{U}^{\rm rr},{\bf Z}/p{\bf Z})$ is $\leqslant {\rm size}(A) \cdot [{\bf F}_q:{\bf Z}/p{\bf Z}]$.

\medskip

Conversely, the group $\overline{U}^{\rm rr}$ is the pro-$p$ Sylow subgroup of the stabilizer of the standard chamber, say $c$, in the geometric completion $\overline{G}^{\rm rr}$ (see end of \ref{ss - KM complete}). 
Let us call $B$ the combinatorial ball of radius 1 around $c$ in the building $X$ acted upon by $\overline{G}^{\rm rr}$.
It follows from the defining relations of minimal Kac-Moody groups, that for each simple root $\alpha$, the root group $U_\alpha$ acts simply transitively on the chambers in $B$ intersecting $c$ along its panel of type $\alpha$.
Similarly, if $\alpha'$ is another simple root, then the chambers not in the standard apartment and intersecting $c$ along its panel of type $\alpha'$ are of the form $u.c$ where 
$u \in U_{-\alpha'}\setminus \{Ê1 \}$.
But since $\{Ê\alpha;-\alpha' \}$ is a prenilpotent pair of real roots whose sum is not a root, it follows again from the defining relations of minimal Kac-Moody groups, that 
$[U_\alpha,U_{-\alpha'}]=\{Ê1\}$, so that $U_\alpha$ fixes pointwise the chambers in $B$ not containing the panel of type $\alpha$ of $c$.
This description of the action of each simple root group on $B$ shows that the quotient of $\overline{U}^{\rm rr}$ corresponding to the restriction of the action of $\overline{U}^{\rm rr}$ on $B$ is an abelian group of exponent $p$ and order $p^{{\rm size}(A) \cdot [{\bf F}_q:{\bf Z}/p{\bf Z}]}$.
Since $\overline{U}^{\rm rr}/\Phi(\overline{U}^{\rm rr})$ is the biggest quotient of $\overline{U}^{\rm rr}$ which is an abelian group of exponent $p$, we conclude that the dimension of ${\rm H}_1(\overline{U}^{\rm rr},{\bf Z}/p{\bf Z})$ is equal to ${\rm size}(A) \cdot [{\bf F}_q:{\bf Z}/p{\bf Z}]$.

\medskip 

Finally, we have ${\rm H}_1(\overline{U}^{{\rm cg}\lambda},{\bf Z}/p{\bf Z}) \simeq ({\bf Z}/p{\bf Z})^{{\rm size}(A) \cdot [{\bf F}_q:{\bf Z}/p{\bf Z}]}$ in view of the interpretation of ${\rm H}_1(\overline{U}^{{\rm cg}\lambda},{\bf Z}/p{\bf Z})$ in terms of topologically generating sets for $\overline{U}^{{\rm cg}\lambda}$ \cite[I.4.2]{SerreCoGa}.
\hfill $\Box$

\medskip

In the classical case, Theorem \ref{th.Frattini} provides the (presumably well-known) fact below.

\medskip

\begin{Cor}
\label{cor.linear}
Let ${\bf G}$ be a split simply connected simple algebraic group over ${\bf F}_q$.
Let $\overline{U}$ be the pro-$p$ Sylow subgroup of the hyperspecial maximal compact subgroup ${\bf G}({\bf F}_q[[t]])$ in ${\bf G}\bigl({\bf F}_q(\!(t)\!)\bigr)$.
Then $\Phi(\overline{U}) = [\overline{U},\overline{U}]$ and ${\rm dim}_{{\bf Z}/p{\bf Z}}Ê\, {\rm H}_1(\overline{U},{\bf Z}/p{\bf Z}) = ({\rm rk}_{{\bf F}_q(\!(t)\!)}({\bf G})+1)\cdot [{\bf F}_q:{\bf Z}/p{\bf Z}]$.
\end{Cor}

\medskip

\begin{Rk}
\label{rk.nonsplit}
This corollary can probably be extended to the non-split isotropic case, but the tools will then be relevant to Bruhat-Tits theory, and not Kac-Moody theory.
\end{Rk}

\medskip

\noindent \emph{\textbf{Proof}}.
This follows from the fact that ${\bf G}\bigl({\bf F}_q(\!(t)\!)\bigr)$ modulo its center is a complete Kac-Moody group over ${\bf F}_q$ with root system the affine root system $\widetilde R$ associated to the finite root system $R$ of ${\bf G}$.
For the root system $\widetilde R$, the corresponding simple root groups are given by the simple root groups in a Chevalley \'epinglage, say $\{Êu_\alpha : \mathbb{G}_a \to U_\alpha \}_{\alpha \in R}$,  together with the group $\{Êu_{-\delta}(rt) : r \in \mathbf{F}_q \}$ where $\delta$ is the longest root in $R$.
\hfill $\Box$
\medskip

\begin{Ex}
\label{Ex.Frattini}
The group ${\rm SL}_m\bigl({\bf F}_q(\!(t)\!) \bigr)$ is a complete Kac-Moody group of affine type.
Any $p$-Sylow subgroup is conjugated to the subgroup $\overline U$ of ${\rm SL}_n({\bf F}_q[[t]])$ containing all the matrices which are unipotent upper-triangular modulo $t$.
We assume that the ground field ${\bf F}_q$ is of degree $r$ over ${\bf Z}/p{\bf Z}$ and pick a corresponding basis $v_1=1$, $v_2$, $\dots$ $v_r$.
Denoting by the $E_{i,j}$ $(1 \leqslant i,j \leqslant m)$ the elements of the canonical basis of the $m \times m$ matrices, we see that the theorem implies that a topologically generating set for $\overline U$ of minimal cardinality is given by the matrices $1 + v_l E_{i,i+1}$ and $1 + v_l E_{m,1}$ $(1 \leqslant l \leqslant r$ and $1 \leqslant i \leqslant m-1)$.
\end{Ex}

\begin{Rk}
Mentioning this linear example is the opportunity to note that there is one, simple but we think convincing, computation to do in order to understand why imaginary roots should finally be taken into account.
Namely, the commutator of elementary unipotent matrices:

\medskip

\centerline{$[\left( \begin{array}{cc}
1 & rt^m  \\
0 & 1  \end{array} \right),
\left( \begin{array}{cc}
1 & 0 \\
st^n & 1  \end{array} \right)]=
\left( \begin{array}{cc}
1 & rt^m  \\
0 & 1  \end{array} \right)
\left( \begin{array}{cc}
1 & 0 \\
st^n & 1  \end{array} \right)
\left( \begin{array}{cc}
1 & -rt^m  \\
0 & 1  \end{array} \right)
\left( \begin{array}{cc}
1 & 0 \\
-st^n & 1  \end{array} \right)$}

\medskip

\noindent is equal to the matrix 

\medskip

\centerline{$
\left( \begin{array}{cc}
1+ rst^{m+n} + r^2s^2t^{2m+2n} & \hfill -r^2st^{2m+n}  \\
\hfill rs^2t^{m+2n} & 1- rst^{m+n}  \end{array} \right)$.}

\medskip 

\noindent The off-diagonal coefficients there correspond to real root groups, but typically the diagonal matrices of the form 
$
\left( \begin{array}{cc}
1 + tP(t)  & 0  \\
0 & \bigl(1 + tP(t)\bigr)^{-1} \end{array} \right)$,
with $P \in \mathbf{F}_q(\!(t)\!)$, correspond to imaginary root groups.
These matrices appear naturally in completions, but not in minimal Kac-Moody groups, even when $P$ is a polynomial.
Using the action of ${\rm SL}_n({\bf F}_q[t])$ on the Bruhat-Tits tree of ${\rm SL}_n\bigl({\bf F}_q(\!(t )\!) \bigr)$, one sees that there is no relation in ${\rm SL}_n({\bf F}_q[t])$ between any element 
of
$\left( \begin{array}{cc}
1 & {\bf F}_q[t]  \\
0 & 1 \end{array} \right)$
and any element of
$\left( \begin{array}{cc}
1 & 0  \\
t{\bf F}_q[t] & 1 \end{array} \right)$ \cite[Proposition 3, p.88]{SerreArbres},
but there are relations in ${\rm SL}_n({\bf F}_q[[t]])$.
Geometrically, these relations appear because ${\rm PSL}_n({\bf F}_q[[t]])$ is the projective limit of finite groups corresponding to the restriction of the action of ${\rm PSL}_n({\bf F}_q[t])$ to finite balls centered at a suitable edge in the tree.
For general complete Kac-Moody groups, imaginary root groups are expected to be very helpful to describe the full pointwise stabilizer of an apartment in the associated building.
\end{Rk}

\medskip

\subsection{Related questions}
\label{ss - fg q}

We think that it would be very interesting to investigate further the higher homology groups ${\rm H}_i(U^{{\rm ma}+},{\bf Z}/p{\bf Z})$. 
In particular, it would be interesting to compute ${\rm H}_2(U^{{\rm ma}+},{\bf Z}/p{\bf Z})$ since it is connected to the size of pro-$p$ presentations of $U^{{\rm ma}+}$ \cite[I.4.3]{SerreCoGa}.

One question in this vein is to understand whether the cohomological finiteness depends on the defining Lie-theoretic data and, maybe, on the ground field.
In the discrete case, there is indeed a dependence on these parameters for cohomology vanishing \cite{DJ}. 
The most obvious point would be to know whether topological finite generation is true without any assumption comparing $p$ or $q$ with respect to $A$.

At last, the proof of the above theorem suggests that the Lie-theoretic nature of $U^{{\rm ma}+}$ will make manageable the computation of some significant sequences of closed subgroups in 
$U^{{\rm ma}+}$ (at least in connection with the relevant infinite root system). 
This may be a hint in the study of linearity of the groups $U^{{\rm ma}+}$ (see Sect. \ref{s - NonLinearity}).

\section{Abstract simplicity for complete Kac-Moody groups}
\label{s - Simplicity}

In this section, we derive abstract simplicity of complete Kac-Moody groups from topological finite generations of the relevant Sylow subgroups.
The question of abstract simplicity had been settled in many interesting cases in \cite{CER}, and is now proved in full generality for geometric completions by T. Marquis \cite[Theorem A]{Marquis} thanks to a very nice argument of dynamical nature.
In our results, given a generalized Cartan matrix $A$, we have to make some assumption comparing the characteristic of the ground field and the coefficients of $A$, but our completions are of arbitrary type.
Our proof is algebraic and uses arguments from Tits systems in the spirit of Bourbaki's simplicity criterion; it owes to \cite{CER} the idea that topological finite generation of pro-$p$ Sylow subgroups implies abstract simplicity.

\subsection{Tits systems}
\label{ss - BN-pair}

Tits system is a structure of combinatorial nature in group theory: it involves subgroups, usually called $B$ and $N$, which are requested to satisfy conditions essentially dealing with double cosets modulo $B$.
More precisely \cite[IV.2.1]{Bourbaki}, a {\it Tits system}~(or $BN$-pair) is a quadruple $(G,B,N,S)$ where $G$ is a group, $B$ and $N$ are two subgroups and $S$ is a subset of $N/(B \cap N)$, satisfying the following axioms.
\begin{itemize}
\item[{\rm (T1)}] The subset $B \cup N$ generates $G$ and $B \cap N$ is normal in $N$.
\item[{\rm (T2)}] The set $S$ generates the group $W = N/(B \cap N)$ and consists of elements of order 2 in $W$.
\item[{\rm (T3)}] We have: $sBw \subset BwB \cup BswB$ for $s \in S$ and $w \in W$.
\item[{\rm (T4)}] For any $s \in S$ we have $sBs \not\subset B$.
\end{itemize}

\medskip

It follows from the axioms that $(W,S)$ is a Coxeter system \cite[IV.2.4, Th\'eor\`eme 2]{Bourbaki}.

\medskip

We will use the following simplicity criterion \cite[IV.2.7, Th\'eor\`eme 5]{Bourbaki}

\medskip

\begin{Thm}
\label{th.BBK}
Let $(G,B,N,S)$ be a Tits system with Weyl group $W = N/(B \cap N)$.
Let $U$ be a subgroup of $B$. 
We set $Z = \bigcap_{g \in G} gBg^{-1}$ and $G^\dagger = \langle gUg^{-1} : g \in G\rangle$.
Assume that the following holds.
\begin{itemize}
\item[{\rm (1)}] We have $U \triangleleft B$ and $B = UT$ where $T = B \cap N$.
\item[{\rm (2)}] For any proper $V \triangleleft U$, we have $[U/V,U/V] \subsetneq U/V$.
\item[{\rm (3)}] We have $G^\dagger=[G^\dagger,G^\dagger]$.
\item[{\rm (4)}] The Coxeter system $(W,S)$ is irreducible.
\end{itemize}

\noindent Then for any subgroup $H$ in $G$ normalized by $G^\dagger$ we have either $H \subseteq Z$ or $G^\dagger \subseteq H$.
\end{Thm}

\medskip

Complete Kac-Moody groups fit in this general picture: indeed, it follows from \cite[Proposition 3.16]{Rousseau} that the algebraic completion $\frak{G}_\mathcal{D}^{\rm ma+}({\bf F}_q)$ admits a refined Tits system in the sense of \cite{KacPet}.
Such a structure implies the existence of a Tits system and the group $U^{{\rm ma}+}$ satisfies condition (1) stated for the group $U$ in the theorem above.
Moreover we know that the associated Weyl group is irreducible if and only if the generalized Cartan matrix we started with is indecomposable.
Concerning the geometric completion, the analogous properties (after replacing $U^{{\rm ma}+}$ by $\overline{U}^{\rm rr}$) had been proved in \cite[Theorem 1.C]{RemRon}.

\medskip

Therefore in what follows, we will essentially have to check conditions (2) and (3) to derive simplicity results from Bourbaki's criterion.
 
\subsection{Abstract simplicity}
\label{ss - simple}

We now turn to the algebraic proof of abstract simplicity for "most "complete Kac-Moody over finite fields.

\medskip

\begin{Thm}
\label{th.simple}
Let $\mathcal{D}$ be a Kac-Moody root datum with generalized Cartan matrix $A$.
Let ${\bf F}_q$ be a finite field of characteristic $p$. 
Let $G$ be the algebraic or geometric completion of the minimal Kac-Moody group $\frak{G}_\mathcal{D}({\bf F}_q)$.
We denote by $G^\dagger$ the group topologically generated by the root groups in $G$ and by $Z'(G)$ the kernel of the action of $G$ on the building given by its standard refined Tits system.
We assume that $A$ is indecomposable, that $q \geqslant 4$ and that $p$ is greater than the absolute value of any off-diagonal coefficient of $A$.
Then any subgroup of $G$ normalized by $G^\dagger$ either is contained in $Z'(G)$ or contains $G^\dagger$.
In particular, the group $G^\dagger/(Z'(G)\cap G^\dagger)$ is abstractly simple.
\end{Thm}

\begin{Rk}
\label{rk.simple.Rousseau}
This statement is the analogue over finite fields of \cite[Theorem 6.19]{Rousseau}. 
\end{Rk}

\medskip

\noindent \emph{\textbf{Proof}}.
We denote by $\overline U$ the standard pro-$p$ Sylow subgroup of $G$.

\medskip

First step: perfectness of $G^\dagger$.

\smallskip

\noindent Since $q \geqslant 4$, the group $G^\dagger$ contains a perfect (often simple modulo center) dense subgroup, 
namely the subgroup of $\frak{G}_\mathcal{D}({\bf F}_q)$ generated by the root subgroups \cite[Lemme 6.14]{Rousseau}, so that $[G^\dagger,G^\dagger]$ is dense in $G^\dagger$.
But we saw in \ref{ss - fg} that $[\overline U, \overline U]=\Phi(\overline U)$ and that the latter group is open in $\overline U$.
This implies that the dense subgroup $[G^\dagger,G^\dagger]$ also contains an open subgroup.
Therefore we have $G^\dagger=[G^\dagger,G^\dagger]$.

\medskip

Second step: non-perfectness of the quotients of $\overline U$.

\smallskip

\noindent Let $V \triangleleft \overline U$ be proper.
We want to see that $[\overline U/V, \overline U/V]Ê\subsetneq \overline U/V$, which amounts to seeing that the inclusion $V \cdot [\overline U, \overline U]Ê\subset \overline U$ is strict.
But if we had the equality, then since $\Phi(\overline U)=[\overline U,\overline U]$, this would mean that $V \cdot \Phi(\overline U) = \overline U$, implying that $V$ generates $\overline U$; this is excluded since $V$ is by assumption a proper subgroup.

\medskip

Conclusion: it remains to apply Bourbaki's simplicity theorem to the refined Tits system of the previous subsection.
\hfill $\Box$

\medskip

\begin{Rk}
\label{history of simplicity}
Topological simplicity is proved for geometric completions in \cite[Theorem 2.A.1]{RemyGAFA}, where the idea of combining arguments from Tits systems and pro-$p$ groups is used but not pushed far enough to obtain abstract simplicity.
The idea to use more subtle pro-$p$ group arguments to obtain abstract simplicity appears in \cite{CER} and the proof here follows the lines of the argument in the latter paper.
At last, the most general result of abstract simplicity is due to T. Marquis.
He uses completely different techniques \cite{Marquis}, namely contraction groups.
Note also that many abstract simplicity results for complete Kac-Moody groups over infinite fields were obtained by R. Moody \cite{Moody} and G. Rousseau \cite[Section 6]{Rousseau}.
\end{Rk}

\medskip

\begin{Ex}
\label{simple.linear}
In the linear case, i.e. when $G={\bf G}\bigl({\bf F}_q(\!(t)\!)\bigr)$ where ${\bf G}$ is a split simply connected simple algebraic group over ${\bf F}_q$ (see Corollary \ref{cor.linear}), the proof above uses a refined Tits system with affine Weyl group and topological finite generation of the pro-$p$ Sylow subgroups of $G$.
Of course, a better proof in this case is the classical one, which uses the usual refined Tits system (with finite Weyl group), in which the group $U$ is the (obviously non-perfect) unipotent radical of a minimal parabolic subgroup.
\end{Ex}

\subsection{Related questions}
\label{ss - simple q}

In the same vein of results, it would be interesting to prove (or disprove) abstract simplicity of geometric completions of groups with mixed ground fields as defined in \cite{RemRon}.
In that case, the tools developed by G.~Rousseau's paper \cite{Rousseau} are no longer available since there is no infinite-dimensional Lie algebra to start with.
Even more exotic examples of groups acting on locally finite Moufang twin buildings, and therefore admitting locally compact geometric completions, were given in \cite{AbrRem}.

\section{Non-linearity statements and conjectures}
\label{s - NonLinearity}

In this section, we disprove some linearities for pro-$p$ Sylow subgroups $\overline U$ of complete Kac-Moody groups.
Actually, we mention some conjectures, and prove only partial results in their direction.
Still, we believe that our questions may be found challenging by experts in pro-$p$ groups.
We mention two approaches. 
The first one uses the notion commensurator for profinite groups \cite{BEW}, combined with some rigidity phenomena implied by R.~Pink's work \cite{Pink}. 
The second one, more shortly developed, is a variant which uses the notion of Golod-Shafarevich groups. 
One obstacle to this second strategy to disprove linearities of our pro-$p$ groups $\overline U$ is the fact, shown below, that in general the groups $\overline U$ differ from the full pro-$p$ completions $\widehat{U}^p$ for which useful results are available \cite{ErshovGoSha}.

\subsection{Further material on profinite groups}
\label{ss - commensurator}

 We begin by recalling few definitions. Let $P$ be a profinite group.
A couple $(L,\eta)$ is called an $\it{envelope}$ of $P$ if it consists of a topological group $L$ and an injective homomorphism $\eta:P\rightarrow L$ such that $\eta(P)$ is open in $L$ and $\eta$ maps $P$ homeomorphically onto $\eta(P)$ \cite[Definition 3.2]{BEW}.
The group $L$ itself is also  referred to as an envelope of $P$ whenever the reference to $\eta$ is inessential.

\medskip 

Consider all the $\it{virtual\ automorphisms}$  of $P$,  i.e. all the topological isomorphisms  from an open subgroup of $P$ to another open subgroup of $P$. We will say
that two virtual automorphisms  of $P$ are equivalent if they coincide on some open subgroup of $P$. Since the intersection of open subgroups is again an open subgroup, one can endow the set of equivalence  classes of virtual automorphisms of $P$ with a natural operation induced by composition. It follows immediately that this is a group which is called the $\it{commensurator}$ of $P$ and is denoted by $\rm{Comm}(P)$.

\medskip

The following rigidity statement is a famous result of R.~Pink \cite[Corollary 0.3]{Pink}.

\medskip

\begin{Thm}
\label{Thm.Pink}
For each $i=1,2$ consider a local field $\mathbf{F}_i$, an absolutely simple, simply connected group $\mathbf G_i$ over $\mathbf{F}_i$ , and an open compact subgroup $U_i\subset \mathbf{G}_i(\mathbf{F}_i)$. Let $f : U_1\rightarrow U_2$ be an isomorphism of topological groups. Then 
there exists a unique isomorphism of algebraic groups $\mathbf{G}_1\rightarrow \mathbf{G}_2$ over a unique 
isomorphism of local fields $\mathbf{F}_1\rightarrow \mathbf{F}_2$ , such that the induced isomorphism 
$\mathbf{G}_1(\mathbf{F}_1)\rightarrow  \mathbf{G}_2(\mathbf{F}_2)$ extends $f$.
\end{Thm}

\medskip 

This theorem is fruitfully combined with the following fundamental result of A.~Borel and J.~Tits \cite[Corollary 8.13]{BorelTits}.

\medskip

\begin{Thm}
\label{Thm.BorelTits}
Let $\mathbf{G}$ be  an absolutely almost simple, simply connected algebraic group defined and isotropic over some infinite field $F$.
Then ${\rm Aut}\bigl( \mathbf{G}(F) \bigr)\cong \bigl({\rm Aut}(\mathbf{G})\bigr)(F)\rtimes {\rm Aut}(F)$ where ${\rm Aut}(\mathbf{G})$  is the group of automorphisms of the algebraic group $\mathbf{G}$.
\end{Thm}

\medskip

Recall also that a pro-$p$ group $V$ is called {\it powerful}~if the quotient $V/\overline{V^p}$ is abelian when $p$ is odd, or if so is $V/\overline{V^4}$ when $p=2$ \cite[Definition 3.1]{DDMS}.
This notion is very useful in the recent reinterpretation of M.~Lazard's work on $p$-adic analytic groups.
This is illustrated for instance by the long list of necessary and sufficient (most of the time, group-theoretic) conditions for a pro-$p$ group to be $p$-adic analytic that can be found in \cite[Interlude A, p.97]{DDMS}.

\subsection{Non-linearities}
\label{ss - non-linear}

In what follows, $\mathcal{D}$ is a Kac-Moody root datum with generalized Cartan matrix $A = [ A_{i,j} ]_{1 \leqslant i,j \leqslant r}$.
Let $\mathbf{F}_q$ be the finite field of characteristic $p$ with $q$ elements.
We denote by $\frak{G}_\mathcal{D}({\bf F}_q)$ the corresponding minimal Kac-Moody group and by $\frak{G}_\mathcal{D}^{\rm ma+}({\bf F}_q)$ its algebraic completion.
We use the notation $U$ for the subgroup of $\frak{G}_\mathcal{D}({\bf F}_q)$ generated by the positive root groups. 

\medskip

Let us start with a first easy non-linearity statement.

\medskip

\begin{Lemma}
\label{lemma.notlinear}
Let $U$ be the group generated by the positive root groups in the minimal Kac-Moody group $\frak{G}_\mathcal{D}(\mathbf{F}_q)$.
Let $\overline U$ be the closure of $U$ in any of the previously considered complete Kac-Moody groups.
Then the group $\overline U$ is not linear over any field of characteristic $\neq p$, hence it is not virtually powerful $($at least when it is finitely generated$)$.
\end{Lemma}

\medskip

\noindent \emph{\textbf{Proof}}.
The first statement follows from the fact that $U$ contains infinite subgroups of exponent $p$ \cite[Theorem 4.6]{RemNotLin}; 
in particular, $\overline U$ cannot be $p$-adic analytic.
When $\overline U$ is finitely generated as a pro-$p$ group, the fact that it is not powerful follows then from \cite[Corollary 8.34]{DDMS}.
\hfill $\Box$

\medskip

\begin{Rk}
\label{rk.linear.pro-p}
Note that the group $\overline U$ is linear when the generalized Cartan matrix $A$ is of affine type since then it is isomorphic to a finite index subgroup of $\mathbf{G}(\mathbf{F}_q[[t]])$ for some suitable linear algebraic group $\mathbf{G}$.
The spherical case is even simpler, since then $\overline U = U$ is a subgroup of a finite group of Lie type.
\end{Rk}

\medskip

The most optimistic conjecture is the following.

\medskip

\begin{Conj}
\label{rk.linear.strong}
{\it With the above notation, assume that $A$ is an indecomposable generalized Cartan matrix whose type is neither spherical nor affine.
Then the group $\overline U$ is not linear over any field.}
\end{Conj}

\medskip

\begin{Rk}
\label{rk.linear.discrete}
Actually, a stronger conjecture shall claim that the countable dense subgroup $U$ is not linear over any field, but it is no longer relevant to pro-$p$ groups: the most serious starting point in this situation seems to be that $U$ acts as a lattice on an explicit ${\rm CAT}(0)$-building for $q >\!\!> 1$.
Note nevertheless that this question is probably much more subtle, since when $A =  \left( \begin{array}{cc}
\hfill 2 & -m  \\
-n & \hfill 2  \end{array} \right)$
with {\it any}~integers $m,n>1$, we have $U \simeq (\mathbf{F}_q[t],+) * (\mathbf{F}_q[t],+)$, which is a linear group (take $m=n=2$).
\end{Rk}

\medskip

As suggested by P.-E. Caprace, the following result can be obtained by using work by R.~Pink \cite{Pink} and Y.~Barnea, M.~Ershov and Th.~Weigel \cite{BEW}.
We keep the previous notation.

\medskip

\begin{Prop}
\label{prop.compact.open}
Assume that $A = [ A_{i,j} ]_{1 \leqslant i,j \leqslant r}$ is an indecomposable generalized Cartan matrix whose type is neither spherical nor affine, with $r \geqslant 3$.
Then the pro-$p$ group $\overline U$ cannot be any compact open subgroup of a simply connected, absolutely almost simple algebraic group defined and isotropic over a non-archimedean local field.
\end{Prop}

\medskip

\begin{Rk}
The assumption on $A$ is made to have $L$ containing an infinite, finitely generated, non-residually finite group, so in view of \cite{CapRemBis}, the conclusion of the Proposition also holds when $A$ is a $2 \times 2$ matrix with one off-diagonal equal to $-1$ and the other of the form $-m$ with $m \geqslant 5$.
\end{Rk}

\medskip

\noindent \emph{\textbf{Proof}}.
We denote by $L$ the ambient complete Kac-Moody group. 
In the terminology of \cite{BEW}, the group $L$ is an envelope of $\overline U$ (\ref{ss - commensurator}); there is a natural group homomorphism $\kappa_L : L \to {\rm Comm}(\overline U)$.
The group $L$ is compactly generated \cite[]{RemyGAFA} and simple \cite[Theorem A]{Marquis}, therefore it follows from \cite[Theorem 4.8]{BEW} that the map $\kappa_L$ is injective.

\medskip

Let $\mathbf{G}$ be a simply connected, absolutely simple algebraic  group defined over some non-archimedean local field $F$.
We want to disprove the possibility that $\overline U$ be isomorphic to a compact open subgroup $\mathbf{G}(F)$.
In view of Lemma \ref{lemma.notlinear} we may -- and shall -- assume that $F \simeq k(\!(t)\!)$ where $k$ is a finite field of characteristic $p$.
Moreover by Theorem~\ref{Thm.Pink} above,  it follows that the commensurator of any compact open subgroup of $\mathbf{G}(F)$ is isomorphic to ${\rm Aut}\bigl( \mathbf{G}(F) \bigr)$. 

\medskip

Hence, let us consider the map $\kappa_L : L \hookrightarrow {\rm Aut}\bigl( \mathbf{G}(F) \bigr)$ in order to obtain a contradiction. 
Since ${\rm Aut}\bigl( \mathbf{G}(F) \bigr)\cong \bigl({\rm Aut}(\mathbf{G})\bigr)(F)\rtimes {\rm Aut}(F)$ by Borel-Tits Theorem \ref{Thm.BorelTits}, we may compose $\kappa_L$ with
the projection $\pi : {\rm Aut}\bigl( \mathbf{G}(F) \bigr) \twoheadrightarrow  {\rm Aut}(F)$ to obtain a map 
from the simple group $L$ to the residually finite group ${\rm Aut}(F)$. It follows immediately that
the map $\pi \circ \varphi$ has to be trivial, meaning that $\kappa_L(L) \subset {\rm Aut}\bigl( \mathbf{G}\bigr)(F)$. 
But ${\rm Aut}\bigl( \mathbf{G}\bigr)(F)$ is a linear group, whereas $L$ contains an infinite finitely generated subgroup (namely the minimal Kac-Moody group generated by the root subgroups) such that the intersection of all its its finite index subgroups still has finite index \cite[Corollary 16]{CapRem}: this is a typical non-linear property which provides the desired contradiction. 

\medskip

Pushing further the use of Pink's result, it might be possible to prove the following.

\medskip

\begin{Conj}
\label{rk.linear.weak}
{\it With the above notation, assume $A$ is an indecomposable generalized Cartan matrix whose type is neither spherical nor affine.
Then for any $m \geqslant 1$ and any local field $F$, there is no injective group homomorphism from $\overline U$ to ${\rm GL}_m(F)$.}
\end{Conj}

\medskip

Arguing by contradiction and starting with an injective group homomorphism $\varphi : \overline U \hookrightarrow {\rm GL}_m(F)$, we know that if $\varphi$ has a relatively compact image then the homomorphism $\varphi$ is automatically continuous \cite[Corollary 1.21]{DDMS}, at least when $\overline U$ is topologically finitely generated.
Therefore the problem splits into two sub-problems: discussing the existence of maps $\varphi$ as above with unbounded image, and showing that $\overline U$ cannot be a compact subgroup in a group ${\rm GL}_m(F)$.
The latter sub-problem is a weak form of Conjecture \ref{rk.linear.strong}, for which \cite[Corollary 0.5]{Pink} may be useful.

\medskip

To finish this subsection, we mention another approach, namely using the fact that free pro-$p$ groups cannot be compact subgroups in linear groups \cite{BarLar}.
The idea would be to combine this fact together with Zelmanov's theorem saying that a Golod-Shafarevich pro-$p$ group contains non-abelian free pro-$p$ groups \cite[Theorem]{ZelmanovHorizons}.
It is already known that the full pro-$p$ completion of some groups $U$ as above has the Golod-Shafarevich property -- surprisingly, some of these discrete groups $U$ enjoy in addition Kazhdan's property (T) \cite{ErshovGoSha}. 
Still, the full pro-$p$ completions seem to be bigger than the closures $\overline U$; unfortunately, we have in fact the following.

\medskip

\begin{Lemma} 
\label{lemma.completion}
There exist complete $($rank $2)$ Kac-Moody groups over $\mathbf{F}_q$ in which the closure $\overline U$ of the subgroup $U$ of the corresponding minimal Kac-Moody group generated by the positive root groups differs from the full pro-$p$ completion $\widehat{U}$ of $U$.
\end{Lemma}

\medskip

\noindent \emph{\textbf{Proof}}.
Suppose  that  the generalized Cartan matrix $A$ involved in the Kac-Moody root datum $\mathcal{D}$ defining the group functor $\frak{G}_\mathcal{D}$ equals
$\left( \begin{array}{cc}
\hfill 2 & -m  \\
-n & \hfill 2  \end{array} \right)$
where $m,n\in\mathbb N$, $m,n>1$ and $p>\max \{Êm;n\}$.
Take $U$,  a subgroup of $\frak{G}_\mathcal{D}({\bf F}_q)$ generated by the positive root groups.
As shown in \cite{Morita}, we have $U=U_1*U_2$ where $U_i$ is an infinite elementary abelian group of exponent $p$ (we have $U_i\cong ({\bf F}_q[t],+)$, $i=1,2$).
Using  \cite[9.1.1]{RZ} we notice that $\widehat{U}=\widehat{U}_1 *_{\hat{p}}
\widehat{U}_2$, a free pro-$p$ product of $\widehat{U}_1$ and $\widehat{U}_2$.
Thus by Lemma 9.1.17 of \cite{RZ}, we have that $\widehat{U}/\Phi(\widehat{U})\cong \widehat{U}_1/\Phi(\widehat{U}_1)\times \widehat{U}_2/\Phi(\widehat{U}_2)\cong\widehat{U}_1\times\widehat{U}_2$. 
However, as we have proved in Theorem~\ref{th.Frattini},  $\overline U$ is a finitely generated pro-$p$ group, which gives an immediate contradiction.
\hfill $\Box$

\subsection{Related questions}
\label{ss - non-linear q}

\medskip

It would be interesting to obtain a description of the quotient groups $\overline U/\overline{(\overline U)^p}$, or $\overline U/\overline{(\overline U)^4}$ when $p=2$.
It follows from \cite{ZelmanovIsrael} that these groups are finite when $\overline U$ is finitely generated, but by Lemma \ref{lemma.notlinear} we know that they are not abelian.
In fact, using Theorem \ref{th.Frattini}, we have  $\Phi(\overline U) = [\overline U, \overline U] \supsetneq \overline{(\overline U)^p}$. 
In some sense, the Lie-theoretic filtrations we have used in this paper are well-adapted to commutators, but are less adapted to $p$-th powers in the groups $\overline U$.

\medskip

A more general question would be to try to characterize pro-$p$ Sylow subgroups arising from {\it non-affine}~Kac-Moody groups with respect to the linear ones.
This could be done by means of the growth of standard sequences of closed subgroups.
This is well-understood for linear pro-$p$ groups in characteristic 0.
The beginning of the investigation of linear groups in positive characteristic (from this viewpoint) appears for instance in \cite{LubSha}.
One idea might be to use the underlying root systems to compute, at least theoretically, some growth rates.

\medskip

\bibliographystyle{amsalpha}
\bibliography{LiePro-p}

\vspace{1cm}

\vspace{0.5cm}
\begin{flushleft} \textit{Bertrand R\'emy} \\
Universit\'e Lyon 1, CNRS - UMR 5208 \\
Institut Camille Jordan \\
43 boulevard du 11 novembre 1918 \\
F-69622 Villeurbanne cedex \\
\vspace{1pt}
remy@math.univ-lyon1.fr
\end{flushleft}

\vspace{0.1cm}
\begin{flushleft}
\textit{Inna Capdeboscq} \\
Mathematics Institute \\
Zeeman Building \\
University of Warwick \\
Coventry CV4 7AL \\

\vspace{1pt}
I.Capdeboscq@warwick.ac.uk
\end{flushleft}

\end{document}